\documentclass[12pt]{amsart}
\usepackage{amsmath, amsthm, amscd, amsfonts}

\newtheorem{theorem}{Theorem}[section]
\newtheorem{lemma}[theorem]{Lemma}

\theoremstyle{definition}
\newtheorem{definition}[theorem]{Definition}

\theoremstyle{remark}
\newtheorem{remark}[theorem]{Remark}
\numberwithin{equation}{section}
\newfont{\kh}{msbm10}

\begin{document}
\title[The Atkinson Theorem in Hilbert $C^*$-modules]{The Atkinson
Theorem in Hilbert $C^*$-modules over $C^*$-Algebras of Compact Operators}
\author{A.Niknam}
\address{Assadollah Niknam, \newline Department of Mathematics,
Ferdowsi University, P. O. Box 1159, Mashhad 91775, Iran;
\newline Centre of Excellence in Analysis on Algebraic
Structures (CEAAS), Ferdowsi Univ., Iran.}
\email{dassamankin@yahoo.co.uk}
\author{K.Sharifi}
\address{Kamran Sharifi, \newline Department of Mathematics,
Ferdowsi University, P. O. Box 1159, Mashhad 91775, Iran}
\email{sharifi.kamran@gmail.com} 
\begin{abstract}
In this paper the concept of unbounded Fredholm operators on
Hilbert $C^*$-modules over an arbitrary $C^*$-algebra is discussed
and the Atkinson theorem is generalized for bounded and unbounded
Feredholm operators on Hilbert $C^*$-modules over $C^*$-algebras
of compact operators. In the framework of Hilbert $C^*$-modules
over $C^*$-algebras of compact operators, the index of an
unbounded Fredholm operator and the index of its bounded
transform are the same.
\end{abstract}
\maketitle
\section{Introduction.}

Hilbert $C^*$-modules are an often used tool in operator theory
and in operator algebras. The theory of Hilbert $C^*$-modules is
very interesting on its own. Interacting with the theory of
operator algebras and including ideas from non-commutative
geometry it progresses and produces results and new problems
attracting attentions (see \cite{GVF}).

A bounded operator $T$ on a Hilbert space $ \mathcal{H}$ is called
Fredholm if there is a bounded operator $S$ on $ \mathcal{H}$ such
that $1-ST$ and $1-TS$ are compact operators on $ \mathcal{H}$.
Atkinson theorem declares that a bounded operator is Fredholm if
and only if its kernel and cokernel are finite dimensional and
its range is closed, see \cite{WEG}, Theorem 14.1.1. The Atkinson
theorem has been extended to bounded operators on standard Hilbert
module $\ell^2( \mathcal{A})$ in section 17.1 of \cite{WEG}. In
the present note we give another generalization of the Atkinson
theorem to Hilbert $C^*$-modules over $C^*$-algebras of compact
operators.

A (left) pre-Hilbert $C^*$-module over a $C^*$-algebra $
\mathcal{A}$ (not necessarily unital) is a left
$\mathcal{A}$-module $E$ equipped with an $\mathcal{A}$-valued
inner product $\langle \cdot , \cdot \rangle : E \times E
\rightarrow \mathcal{A}$, which is $ \mathcal{A}$-linear in the
first variable and has the properties:$$ \langle x,y
\rangle=\langle y,x \rangle ^{*}, \ \langle x,x \rangle \geq 0 ~
\textrm{with equality if and only if}~ x=0$$

A pre-Hilbert $ \mathcal{A}$-module $E$ is called a Hilbert $
\mathcal{A}$-module if it is complete with respect to the norm $\|
x \|=\|\langle x,x\rangle\| ^{1/2}$. The basic theory of Hilbert
$C^*$-modules can be found in \cite{LAN}, \cite{WEG}.

We denote by $B(E)$ the $C^*$-algebra of all adjointable operators
on $E$ (i.e. of all maps $T :E\rightarrow E$ such that there
exists $T^* :E\rightarrow E$ with the property $\langle Tx,y
\rangle =\langle x,T^*y \rangle$, for all $ x,y\in E$). It is
well known that each adjointable operator is necessarily bounded
$\mathcal{A}$-linear in the sense $T(ax) = aT(x)$, for all $a\in
\mathcal{A}$, $x\in E$. In general, bounded $ \mathcal{A}
$-linear operator may fail to possess an adjoint (cf. \cite{LAN}).
However, if $E$ is a Hilbert $C^*$-module over the $C^*$-algebra
$ \mathcal{K} $ of all compact operators on a Hilbert space, then
it is known that each bounded $ \mathcal{K}$-linear operator on
$E$ is necessarily adjointable (see for example \cite{B-G},
Remark 5).

Given elements $x,y\in E$ we define $\Theta _{x,y} :E\rightarrow
E$ by $\Theta_{x,y}(z) = \langle z,x \rangle y $ for each $z \in
E$, then $\Theta_{x,y} \in B(E)$, with $(\Theta_{x,y}
)^*=\Theta_{y,x}$. The closure of the span of $\{ \Theta_{x,y}:
x,y \in E \}$ in $B(E)$ is denoted by $K(E)$, and elements from
this set will be called $\mathcal{A}$-compact operators.

In various contexts where the Hilbert $C^*$-modules arise, one
also needs to study "unbounded adjointable operators" or what are
now know as regular operators. These were first introduced by
Baaj and Julg in \cite{B-J} where they gave a nice construction
of Kasparov bimodules in $KK$-theory using regular operators,
Lance gave a brief indication in his book \cite{LAN} about
Hilbert modules and regular operators on them. Let us quickly
recall the definition of a regular operator. An operator $T$ from
a Hilbert $ \mathcal{A}$-module $E$ to $E$ is said to be regular
if

(i) $T$ is closed and densely defined,

(ii) its adjoint $T^*$ is also densely defined, and

(iii) range of $1+T^*T$ is dense in $E$.

Note that as we set $ \mathcal{A}= \mathbb{C} $ i.e. if we take
$E$ to be a Hilbert space, then this is exactly the definition of
a closed operator, except that in that case, both the second and
third condition follow from the first one. Before starting the
other parts let us fix the rest of our notation.

Throughout the paper, $ \mathcal{K}=\mathcal{K( \mathcal{H})}$
will be the $C^*$-algebra of all compact operators on a Hilbert
space $ \mathcal{H}$ and $ \mathcal{A}$ will be an arbitrary
$C^*$-algebra (not necessarily unital). In this paper we will
deal with bounded and unbounded operators at the same time. To
simplify notation we will, as a general rule, denote bounded
operators by capital letters and unbounded operators by small
letters, also we will use $Dom(.)$ for the domain of unbounded
operators. We use $Ker$ and $Ran$ for the
kernel and the range of operators respectively.\\
\section{Preliminaries}

In this section we would like to recall some definitions and
present a few simple facts about unbounded regular operators and
their bounded transforms. We give a definition of regular
Fredholm operators on Hilbert $C^*$-modules and then we state
that a regular operator is Fredholm if and only if its bounded
transform is a Fredholm operator.

Let $E$ be a Hilbert $ \mathcal{A}$-module, consider a closed $
\mathcal{A}$-linear operator $t :Dom(t) \subseteq E \rightarrow
E$ where $Dom(t)$ is a dense submodule of $E$. We define $$
Dom(t^*)= \{y \in E : \exists \ y^{'} \in E \ \textrm{s.t.} \
\langle tx,y \rangle = \langle x , y^{'}\rangle \ \textrm{for all
}~ x \in Dom(t) \} .$$This is the domain of a closed $
\mathcal{A}$-linear operator $ t^* :Dom(t^*) \subseteq E
\rightarrow E $ uniquely determined by $ \langle tx,y \rangle =
\langle x,t^*y \rangle$ for all $x \in Dom(t), \ y \in Dom(t^*) .
$
\begin{definition} An operator $t$ as above is called regular if
$t^*$ is densely defined and $1+t^*t$ has dense range. The set of
all regular operators on $E$ is denoted by $R(E)$.
\end{definition}
If t is a regular operator so are $t^*$ and $t^*t$, moreover
$t^{**}=t$ and $t^*t$ is self adjoint (cf. \cite{LAN}, Corollaries
9.4, 9.6 and Proposition 9.9), and also we can define $$ F_{t}=
t(1+t^*t) ^{-1/2}$$ $$ Q_{t}=(1+t^*t) ^{-1/2}$$ then $F_{t}$,
$Q_{t}$ and $tQ_{t}^{2}$ are in $B(E)$ and $ Ran\, Q_{t}=Dom(t)$
(cf. \cite{LAN}, chapter 9).

The map $t\rightarrow F_{t}$ defines a bijection (cf. \cite{LAN},
Theorem 10.4)
$$R(E) \rightarrow \{F\in B(E):\|F\|\leq 1 \ \textrm{and} \
Ran(1-F^*F) \ \textrm{is dense in}~ E \}. $$ This map is
adjoint-preserving i.e. $F_{t}^*=F_{t^*}$ and
$F_{t}=tQ_{t}=t(1+t^*t) ^{-1/2}$ is called the bounded transform
of a regular operator $t$. Moreover we have $
Q_{t}=(1-F_{t}^*F_{t})
^{1/2}$ and $t=F_{t}(1-F_{t}^*F_{t}) ^{-1/2},\|F_{t}\|\leq 1  \ .$\\

Recall that the composition of two regular operator $t,s\in R(E)$
is the unbounded operator $ts$ with $Dom(ts)=\{x\in Dom(s) : \
sx\in Dom(t) \}$ given by $(ts)x=t(sx)$ for all  $x\in Dom(ts)$.
Note that, in general, the composition of regular operators will
not be a regular operator.

Recall that a bounded operator $T\in B(E)$ is said to be Fredholm
(or $ \mathcal{A}$-Fredholm) if it has a pseudo left as well as
pseudo right inverse i.e. there are $S_{1},S_{2}\in B(E)$ such
that $S_{1}T=1 \ \ mod \ K(E)$ and $TS_{2}=1 \ \ mod \ K(E)$.
This is equivalent to say that there exist $S\in B(E)$ such that
$ST=TS=1 \ \ mod \ K(E)$. For more details about bounded Fredholm
operators on Hilbert $C^*$-modules and their applications one can
see \cite{GVF}, \cite{WEG}. The theory of unbounded Fredholm
operators on Hilbert spaces and on standard Hilbert $
\mathcal{A}$-module $ \ell^{2}( \mathcal{A})$ are discussed in
\cite{BLP} and \cite{JOA} respectively. These motivate us to study
such operators on general Hilbert $C^*$-modules as follow:
\begin{definition} Let $t$ be a regular operator on a Hilbert $
\mathcal{A}$-module $E$. An adjointable bounded operator $G\in
B(E)$ is called a pseudo left inverse of t if $Gt$ is closable
and its closure $ \overline{Gt}$ satisfies $ \overline{Gt}\in
B(E)$ and $ \overline{Gt}=1 \ \ mod \ K(E)$. Analogously $G$ is
called a pseudo right inverse if $tG$ is closable and its closure
$ \overline{tG}$ satisfies $ \overline{tG}\in B(E)$ and $
\overline{tG}=1 \ \ mod \ K(E)$. The regular operator t is called
Fredholm (or $ \mathcal{A}$-Fredholm), if it has a pseudo left as
well as a pseudo right inverse.
\end{definition}
When we are dealing with Fredholm operators, a useful connection
between unbounded operators and their bounded transforms (on $
\ell^{2}( \mathcal{A})$) has been given in \cite{JOA}, Lemma 2.2.
Fortunately its proof can be repeated word by word to get the
following result.
\begin{theorem}Let $E$ be a Hilbert $  \mathcal{A}$-module and $t$
be a regular operator on $E$. Then $t$ is Fredholm if and only if
$F_{t}$ is.
\end{theorem}
\begin{proof}Reffer to \cite{JOA}, Lemma 2.2.
\end{proof}

\section {Bounded Fredholm operators}

In this section we introduce a concept of an orthonormal basis
for Hilbert $C^*$-modules and then we briefly discuss bounded
Fredholm operators in Hilbert $ \mathcal{K}$-modules. For this aim
we borrow some definitions from \cite {B-G}.

Let $E$ be a Hilbert $ \mathcal{A}$-module, a system
$(x_{\lambda}), \lambda \in \Lambda $ in $E$ is orthonormal if
each $ x_{\lambda}$ is a basic vector (i.e. $e=\langle
x_{\lambda},x_{\lambda} \rangle$ is minimal projection in $
\mathcal{A}$) and $ \langle x_{\lambda},x_{\mu} \rangle=0$ for
all $\lambda \neq \mu$. An orthonormal system  $( x_{\lambda})$ in
$E$ is said to be an orthonormal basis for $E$ if it generates a
dense submodule of $E$. Immediately, the previous definition
implies that if $x\in E$ satisfies $\langle x,x\rangle=e$ for some
projection (not necessarily minimal) $e\in \mathcal{A} $ then
$\langle ex-x,ex-x\rangle =0$ so $ex=x$. In particular, the same
is true for all basic vectors in $E$.
\begin{theorem}Let $E$ be a Hilbert $ \mathcal{K}$-module, then
there exists an orthonormal basis for $E$.
\end{theorem}
\begin{proof}Reffer to \cite{B-G}, Theorem 4.
\end{proof}
\begin{definition}Let $E$ be a Hilbert $ \mathcal{K}$-module. The
orthonormal dimension of $E$ (denoted by $dim_{ \mathcal{K}} E$)
is defined as the cardinal number of any of its orthonormal bases.
\end{definition}

Note that every two orthonormal bases of $E$ have the same
cardinal number (see \cite{B-G}).

Now, let us fix a minimal projection $e_{0}\in \mathcal{K}$ and
denote by $E_{e_{0}}=e_{0}E=\{e_{0}x \ : \ x\in E \}$ then
$E_{e_{0}}$ is a Hilbert space with respect to the inner product
$(.,.)=tr(\langle.,.\rangle)$, which is introduced in \cite{B-G},
Remark $4$, moreover $dim_{ \mathcal{K}} E=dim \,E_{e_{0}}$.

Suppose $B(E_{e_{0}})$ and $K(E_{e_{0}})$ are the $C^*$-algebras
of all bounded linear operators and compact operators on Hilbert
space $E_{e_{0}}$, respectively, then $B(E)$ and $K(E)$ can be
described by the following theorem.

\begin{theorem}Let $E$ be a Hilbert $ \mathcal{K}$-module and let $e_{0}$ be an
arbitrary minimal projection in $ \mathcal{K}$. Then the map $\Psi
: $B(E)$\rightarrow B(E_{e_{0}}) , \ \Psi(T)= T| _{ E_{e_{0}} } $
is a $*$-isomorphism of $C^*$-algebras. Moreover $T$ is a compact
operator on $E$ if and only if $\Psi(T)= T| _{E_{e_{0}} } $ is a
compact operator on the Hilbert space $E_{e_{0}}$.
\end{theorem}
\begin{proof}Reffer to \cite{B-G}, Theorems 5, 6.
\end{proof}

\begin{remark}Let $E$ be a Hilbert $ \mathcal{K}$-module and $X$
be a closed submodule in $E$. It is well known that $X$ is
orthogonally complemented in $E$, i.e. $E=X \oplus X ^{\perp}$
(cf. \cite {MAG}).
\end{remark}
\begin{remark}Let $e$ be an arbitrary minimal projection in $
\mathcal{K}$ and suppose $\Psi$ has the same meaning as in
Theorem\,3.3 then $ Ran \, \Psi(T) =e \, Ran\,T$ and $Ker\,
\Psi(T) =e\, Ker \,T $.
\end{remark}
\begin{remark}Let $E$ be a Hilbert $ \mathcal{K}$-module and $T\in
B(E)$, like in the general theory of Banach spaces, one can easily
see that every bounded below operator $T$ on $(Ker \,T)^{\perp} $
has closed range (cf. \cite{MUR}, page 21). It is easily checked
that $T :E\rightarrow E$ is not bounded below on $(Ker
\,T)^{\perp} $ if and only if there is a sequence of unit
elements $x_{n}$ in $(Ker \,T)^{\perp} $ such that $\lim_{n\to
\infty} \,T x_{n}=0 $. This fact will be used in the following
theorem.
\end{remark}
\begin{theorem}Let $E$ be a Hilbert $ \mathcal{K}$-module and
$T\in B(E)$. Then $T$ is Fredholm if and only if the range of $T$
is a closed submodule and both $dim_{ \mathcal{K}} Ker\, T$,
$dim_{ \mathcal{K}} Ker\, T^*$ are finite.
\end{theorem}
\begin{proof}Let $e$ be a minimal projection in $ \mathcal{K}$ and
let \, $\Psi : B(E) \rightarrow B(E_{e})$ be the isomorphism from
Theorem 3.3. Suppose $T$ is a Fredholm operator on $E$ then there
exist an operator $S\in B(E)$ and two compact operator $K_{1}\, ,
K_{2} \in K(E)$ such that $TS-1=K_{1}$ and $ST-1=K_{2}$. Since
$\Psi( K(E))=K(E_{e})$ the operator $\Psi(T)$ is Fredholm on the
Hilbert space $E_{e}$. In particular, $\Psi(T)$ has a closed
range and $Ker\,\Psi(T)$ and $Ker\,\Psi(T^*)$ are finite
dimensional, so by applying Remark 3.5 to the Hilbert modules
$Ker\, T$ and $Ker\ T^*$, respectively, we get $dim_{
\mathcal{K}} Ker\, T$, $dim_{ \mathcal{K}} Ker\, T^*$\,$ <
\infty$.

To prove that range of $T$ is closed it is enough to show that
$T|_{(Ker\,T) ^{\perp} }$ is bounded from below.

Suppose, in the contrary, $T|_{(Ker\,T) ^{\perp}}$ is not bounded
below. Then there  exist a sequence of unit elements $(x_{n}) \in
(Ker\,T) ^{\perp}$ such that $\lim_{n\to \infty} \,T x_{n}=0 $.
This implies $\lim_{n\to \infty} \,T (ex_{n})=e\, \lim_{n\to
\infty} \,T x_{n}=0$, and since $(ex_{n})\in (Ker\, \Psi(T))
^{\perp}$ and the range of $\Psi(T)$ is closed, we obtain a
contradiction.

Conversely, suppose $T \in B(E)$ has closed range and both $dim_{
\mathcal{K}} Ker\, T$ and $dim_{ \mathcal{K}} Ker\, T^*$ are
finite, then Remark 3.5 implies that the range of $\Psi (T)$ is
closed and both $Ker \, \Psi(T)$ and $Ker \, \Psi(T^*)$ are finite
dimensional. Therefore $ \Psi(T)$ is a Fredholm operator on the
Hilbert space $E_{e}$, that is there exists an operator $S \in
B(E_{e})$ such that $ \Psi(T)S-1, S\Psi(T)-1 \in K(E_{e})$.
Utilizing the mapping $ \Psi ^{-1}$, we have $T \Psi^{-1}(S)-1,
\Psi^{-1}(S)T-1 \in K(E)$ and $T$ is therefore a Fredholm operator
on $E$.
\end{proof}

\begin{definition}Let $E$ be a Hilbert $ \mathcal{K}$-module and $T\in
B(E)$ be a Fredholm operator. The Fredholm index of $T$ is an
integer defined by $$ ind\,T=dim_{ \mathcal{K}} Ker\, T - dim_{
\mathcal{K}} Ker\, T^*$$
\end{definition}

The preceding discussion shows how the theory of Fredholm
operators on Hilbert $ \mathcal{K}$-modules is reduced to the
classical theory of Fredholm operators on Hilbert spaces. In
fact, if $E$ is a Hilbert $ \mathcal{K}$-module all properties of
an operator $T\in B(E)$ can be deduced by simple procedure: first,
consider a Hilbert space operator $\Psi(T)$ (by using Theorem
3.3) and then lift the relevant information back to $B(E)$. This
enable us to conclude some results, for example, if $T \in B(E)$
then $T$ is Fredholm operator with $ind \,T=0$ if and only if $T$
is a compact perturbation of an invertible operator.

\section {Unbounded Fredholm operators}

Let $E$ be a Hilbert $ \mathcal{K}$-module. We recall that a
densely defined closed operator  $t :Dom(t) \subseteq E
\rightarrow E$ is said to be regular if $t^*$ is densely defined
and $1+t^*t$ has dense range, moreover $$ F_{t}= tQ_{t}=t(1+t^*t)
^{-1/2} \in B(E) \,, \ F_{t}^*=F_{t^*},
$$ $$ Q_{t}=(1+t^*t) ^{-1/2} \in B(E) \,, \ Ran\, Q_{t}=Dom(t).$$

In this section we are going to prove the Atkinson theorem for
unbounded regular Fredholm operators.

\begin{lemma} Let $E$ be a Hilbert $ \mathcal{K}$-module and $t :Dom(t) \subseteq
E \rightarrow E$ be a regular operator, then $Ker\,t=\{x\in
Dom(t) \, :tx=0 \} \ \,  and \ \, Ker\,t^*=\{x\in Dom(t^*) \,
:t^*x=0 \}$ are closed submodules of $E$.
\end{lemma}
\begin{proof} Let $(x_{n})$ be a sequence in $Ker\,t$ and
$x_{n} \rightarrow x$, then $t(x_{n})=0$, $\forall\, n\in
\mathbb{N}$. Therefore $x_{n}\rightarrow x$ and
$t(x_{n})\rightarrow 0$. It follows, by the closedness of $t$,
that $x\in Dom(t)$ and $tx=0$, that is $x\in Ker\,t$ and so
$Ker\,t$ is closed. Since $t$ is regular, so is $t^*$ and
similarly $Ker\,t^*$ is a closed submodule of $E$.
\end{proof}
\begin{lemma}Let $t$ be a regular operator on a Hilbert $
\mathcal{K}$-module $E$, then

(i) $Ran \, t=Ran \, F_{t}$ and $ Ran \, t^*=Ran \, F_{t^*},$

(ii)\,$Ker\,t^*=(Ran \ t)^{\bot}$ and $Ker\, t=(Ran \
t^*)^{\bot},$

(iii)\,$Ker\, t=Ker\, F_{t}$ and $Ker\, t^*=Ker\, F_{t^*}$.
\end{lemma}
\begin{proof} (i) Recall that $F_{t}=tQ_{t}$ and $ Ran\, Q_{t}=Dom(t)$
then $Ran \, t=Ran \, F_{t}$. Since $t$ is regular and so is
$t^*,$ thus $ Ran \, t^*=Ran \, F_{t^*}.$

(ii) We notice that $y\in Ker\, t^*$ if and only if $ \langle
tx,y \rangle=\langle x,0 \rangle=0$ for all $x\in Dom(t)$, or if
and only if $y\in (Ran\, t) ^{\bot}$. Thus we have $Ker\,t^*=(Ran
\ t)^{\bot}$. The second equality follows from the first equality
and Corollary 9.4 of \cite{LAN}.

(iii) By Theorem 15.3.5 of \cite{WEG}, we have $Ker\,
F_{t^*}=(Ran\, F_{t}) ^{\perp}$ and therefore
$$Ker\, F_{t^*}=(Ran\, F_{t}) ^{\bot}=(Ran\, t)
^{\bot}=Ker\,t^*.$$ Similarly we have $Ker\, F_{t}=Ker\, t.$
\end{proof}
Now we are ready to prove the main theorem of this paper:
\begin{theorem}Let $E$ be a Hilbert $ \mathcal{K}$-module and $t$
be a regular operator on $E$. Then $t$ is Fredholm if and only if
the range of $t$ is a closed submodule of $E$ and both $dim_{
\mathcal{K}} Ker\, t$, $dim_{ \mathcal{K}} Ker\, t^*$ are finite.
\end{theorem}

\begin{proof} By Lemma 4.2, we have $dim_{ \mathcal{K}} Ker\,
t= dim_{ \mathcal{K}} Ker\, F_{t}$,  $dim_{ \mathcal{K}} Ker\,
t^*= dim_{ \mathcal{K}} Ker\,F_{t^*}$ and since
$Ran\,t=Ran\,F_{t}$, the statement is deduced from Theorem 2.3,
and Theorem 3.7.
\end{proof}

Let $t$ be a regular Fredholm operator on Hilbert $
\mathcal{K}$-module $E$ then we can define an index of $t$
formally, that is, we can define: $$ ind \,t=dim_{ \mathcal{K}}
Ker\, t - dim_{ \mathcal{K}} Ker\, t^*,$$ and since $dim_{
\mathcal{K}} Ker\, t= dim_{ \mathcal{K}} Ker\, F_{t}$, $dim_{
\mathcal{K}} Ker\, t^*= dim_{ \mathcal{K}} Ker\,F_{t^*}$ we have
$ind \,t=ind \,F_{t}$.

{\bf Acknowledgement}: The authors would like to thank the
referees for their valuable comments. They also wish to thank
Professor M. S. Moslehian who suggested some useful comments.



\begin{thebibliography}{99}
\bibitem {B-J} S. Baaj and P. Julg, Th\'{e}orie bivariante de
Kasparov et op\'{e}rateurs non born\'{e}s dans les $C^*$-modules
hilbertiens, \textit{C. R. Acad. Sc. Paris S\'{e}r. I Math.}
\textbf{ 296} (1983), no. 21,  875-878.

\bibitem {B-G} D. Baki\'{c} and B. Gulja\v{s}, Hilbert $C^*$-modules over
$C^*$-algebras of compact operators, \textit{Acta Sci. Math.}
(Szeged) \textbf{68} (2002), no. 1-2, 249-269.

\bibitem {BLP} Boss-Bovenbek, M. Lesch and J. Phillips, Unbounded
Fredholm operators and spectral flow, \textit{Canada. J. Math.}
\textbf{57} (2005), no. 2, 225-250.

\bibitem {GVF} J. M. Garcia-Bond\'{\i}a, J. C. V\'{a}rilly and H.
Figueroa, \textit{Elements of non-Commutative geometry},
Birkh\"{a}user 2000.
\bibitem {JOA} M. Joachim, Unbounded Fredholm operators and
K-theory, \textit{High-dimensional manifold topology}, 177-199,
\textit{World Sci. Publishing, River Edge, NJ}, (2003).

\bibitem {LAN} E. C. Lance, \textit{Hilbert $C^*$-Modules}, LMS
Lecture Note Series \textbf{210}, Cambridge Univ. Press, 1995.

\bibitem {MAG} M. Magajna, Hilbert $C^*$-modules in which all
closed submodules are complemented, \textit{Proc. Amer. Math.
Soc.}, \textbf{125} (1997), no. 3, 849-852.

\bibitem {MUR} G. J. Murphy, \textit{$C^*$-algebras and Operator Theory},
Academic Press, 1990.

\bibitem {WEG} N. E. Wegge-Olsen, \textit{K-theory and $C^*$-algebras: a
Frienly Approach}, Oxford University Press, Oxford, England, 1993.
\end{thebibliography}
\end{document}